\documentclass{amsart}
\usepackage{amsfonts,enumerate,amssymb,amsmath,amsthm}
\newtheorem{conjecture}{Conjecture}

\newtheorem{question}{Question}

\newtheorem{theorem}{Theorem}

\newcommand{\threebar}[1]{{\left\vert\kern-0.25ex\left\vert\kern-0.25ex\left\vert #1 
    \right\vert\kern-0.25ex\right\vert\kern-0.25ex\right\vert}}

\usepackage{xcolor}

\newcommand{\R}{\mathbb{R}}
\newcommand{\N}{\mathbb{N}}

\begin{document}
\title[Marginally unstable switched systems with irregular growth]{Marginally unstable discrete-time linear switched systems with highly irregular trajectory growth}
\author{Ian D. Morris}
\email{i.morris@qmul.ac.uk}
\address{School of Mathematical Sciences, Queen Mary, University of London, Mile End Road, London E1 4NS, U.K.}
\maketitle

\begin{abstract}
We investigate the uniform stability properties of discrete-time linear switched systems subject to arbitrary switching, focusing on the ``marginally unstable'' regime in which the system is not Lyapunov stable but in which trajectories cannot escape to infinity at exponential speed. For a discrete linear system of this type \emph{without} switching the fastest-growing trajectory  must grow as an exact polynomial function of time, and a significant body of prior research has focused on investigating how far this intuitive picture can be extended from systems without switching to cases where switching is present. In this note we give an example of a family of discrete linear switched systems in three dimensions, with two switching states, for which this intuition fails badly: for a generic member of this family the maximal rate of uniform growth of escaping trajectories  can be made arbitrarily slow along one subsequence of times and yet also faster than any prescribed slower-than-linear function along a complementary subsequence of times. Using this construction we give new counterexamples to a conjecture of Chitour, Mason and Sigalotti and obtain a negative answer to a related question of Jungers, Protasov and Blondel. Our examples have the additional feature that marginal stability and marginal instability are densely intermingled in the same parameter space.
 \end{abstract}

\section{Introduction and background}
We consider the discrete-time linear switched system whose trajectories are given by
\begin{equation}\label{eq:fundy}x(t+1)=A_{\sigma(t+1)} x(t),\qquad x(0)=v_0\end{equation}
for all integers $t \geq 0$, where each $x(t)$ is a vector in $\R^d$, each $A_i$ is a matrix belonging to some specified subset $\mathsf{A}=\{A_0,\ldots,A_{N-1}\}$ of the vector space $M_d(\R)$ of all real $d \times d$ matrices, $v_0$ is some initial vector in $\R^d$, and $\sigma \colon \N\to \{0,\ldots,N-1\}$ is some function called the \emph{switching law}. We refer to the system \eqref{eq:fundy} in general as the linear switched system defined by $\mathsf{A}$.   

This note is concerned with the worst-case behaviour of the value of the norm $\|x(t)\|$ when the switching law and (normalised) initial vector $v_0$ are chosen arbitrarily. (Here and throughout this note $\|\cdot\|$ will always denote either the Euclidean norm on $\R^d$ or the corresponding operator norm on matrices as appropriate, and for the remainder of the note $v_0$ will always be assumed to be a unit vector.) The growth behaviour of trajectories of discrete-time linear switched systems at exponential scales, as expressed via the concepts of the \emph{joint spectral radius} and \emph{lower spectral radius}, has been the subject of extensive study in the control theory literature and elsewhere (see for example \cite{AhJuPaRo14,Ba88,BoMo15,BrTo79,Gu95,Ju09} and references therein). In this note we will be interested in understanding the worst-case growth of $\|x(t)\|$ at a finer level of detail than its rate of exponential growth. In this we continue the theme of works such as \cite{ChMaSi12,GuZe01,JuPrBl08,PrJu15,Su08}. 

Let us formulate our problem precisely. Given a finite set of matrices $\mathsf{A}=\{A_0,\ldots,A_{N-1}\}\subset M_d(\R)$, how far from the origin can the trajectory $x(t)$ reach at a given time $t\geq 0$, assuming that the initial unit vector $v_0$ and switching sequence $\sigma \colon \N \to \{0,\ldots,N-1\}$ are chosen arbitrarily? In other words, given the set $\mathsf{A}=\{A_0,\ldots,A_{N-1}\}\subset M_d(\R)$, what is the behaviour of the sequence
\begin{equation}\label{eq:banan}t \mapsto \sup_{\sigma \colon \N  \to \{0,\ldots,N-1\}}  \sup_{\substack{v_0 \in \R^d \\ \|v_0\|=1}}  \|x(t)\|?\end{equation}
It is not difficult to see that since we always have
\[x(t)=A_{\sigma(t)}\cdots A_{\sigma(1)}x(0)=A_{\sigma(t)}\cdots A_{\sigma(1)}v_0\]
the sequence \eqref{eq:banan} may be rewritten as the more succinct expression
\begin{equation}\label{eq:that-sequence}t\mapsto\max_{i_1,\ldots,i_t \in \{0,\ldots,N-1\}} \left\|A_{i_t}\cdots A_{i_1}\right\|\end{equation}
and we will often prefer this formulation in the sequel. 

As remarked above, the first natural question concerning the sequence \eqref{eq:that-sequence} is whether or not it grows exponentially. The \emph{joint spectral radius} of the set $\mathsf{A}$ is defined to be the quantity
\[\varrho(\mathsf{A}):=\lim_{t \to \infty} \max_{i_1,\ldots,i_t \in\{0,\ldots,N-1\}} \|A_{i_t}\cdots A_{i_1}\|^{\frac{1}{t}}.\]
This limit always exists and its value is unaffected by a change of norm or a change of basis (for a proof of this and the following statements see \cite{Ju09}). Except in the degenerate case where the elements of $\mathsf{A}$ are all upper triangular matrices with zero diagonal (or where this situation can be achieved by a change of basis for $\R^d$) the joint spectral radius is nonzero. In this degenerate case the sequence \eqref{eq:that-sequence} is identically zero for all $t \geq d$ and for the remainder of this note we will assume $\varrho(\mathsf{A})$ to be nonzero. By replacing the set $\mathsf{A}=\{A_0,\ldots,A_{N-1}\}$ with the set $\mathsf{A}':=\{\varrho(\mathsf{A})^{-1}A_0,\ldots,\varrho(\mathsf{A})^{-1}A_{N-1}\}$ if necessary, when studying the sequence \eqref{eq:that-sequence} we may assume without loss of generality that $\varrho(\mathsf{A})=1$, since any uniform bound $\|x(t)\|\leq P(t)$ on the trajectories of $\mathsf{A}'$ at some time $t \geq 0$ is clearly equivalent to the bound $\|x(t)\|\leq \varrho(\mathsf{A})^tP(t)$ on the trajectories of $\mathsf{A}$. For the remainder of this article we therefore focus on the situation where $\varrho(\mathsf{A})=1$.

When $\varrho(\mathsf{A})=1$ the literature distinguishes two cases: the \emph{marginally stable} case in which there is a uniform bound $\|x(t)\|\leq C$ for all trajectories (where we recall that $x(0)$ is assumed to be a unit vector) and the \emph{marginally unstable} case in which no such uniform upper bound holds\footnote{In this terminology we follow \cite{ChMaSi12,Pr16,PrJu15,Su08} but the term ``marginally unstable'' is in some cases used differently: see for example \cite{ShWiMaWuKi07}.}. In the marginally stable case (which by the homogeneity of the linear maps $A_i$ is precisely the case where the system is Lyapunov stable) it is possible to define a norm on $\R^d$ with respect to which \eqref{eq:that-sequence} is identically equal to $1$ (see \cite{Ju09}) and this implies that in this case \eqref{eq:that-sequence} as measured by the Euclidean norm must be bounded away from both zero and infinity. Consequently, by normalising the joint spectral radius of $\mathsf{A}$ to $1$ and disregarding the already-solved marginally stable case, the analysis of the asymptotic behaviour of \eqref{eq:that-sequence} in general reduces to the question: how can \eqref{eq:that-sequence} behave asymptotically in the marginally unstable case?

When $\varrho(\mathsf{A})=1$ and $\mathsf{A}$ consists of a single matrix $A_0$ (a situation which amounts to the absence of any switching) the behaviour of \eqref{eq:that-sequence} is straightforward. Since we simply have $A_{i_t}\cdots A_{i_1}=A_0^t$ for every $t \geq 1$, either \eqref{eq:that-sequence} grows asymptotically as a polynomial $t \mapsto t^{k}$ for some integer $k<d$ or it is bounded, depending on whether or not there exists an eigenvalue of $A_0$ on the unit circle for which the associated Jordan block is nontrivial. It is therefore natural to ask whether the asymptotic behaviour of \eqref{eq:that-sequence} is also constrained to be either polynomial or bounded in the case $N \geq 2$ in which switching can meaningfully be said to be present. Given $\mathsf{A}=\{A_0,\ldots,A_{N-1}\}\subset M_d(\R)$, let $L_{\mathsf{A}} \geq 1$ denote the largest integer $\ell$ such that in some basis for $\R^d$ we may write
\begin{equation}\label{eq:bananas}A_k=\begin{pmatrix} A^{(1,1)}_k & A^{(1,2)}_k& \cdots & A^{(1,\ell)}_k \\
0&A^{(2,2)}_k & \cdots & A^{(2,\ell)}_k \\
\vdots&\vdots&\ddots&\vdots\\
0&0&\cdots&A^{(\ell,\ell)}_k\end{pmatrix}
\end{equation}
for every $A_k \in \mathsf{A}$, where the matrices along the main diagonal are square and where the dimensions of each block matrix $A^{(i,j)}_k$ depend only on the position $(i,j)$ and not on the particular matrix $A_k \in \mathsf{A}$. Clearly a largest such integer always exists and lies in the range $1,\ldots,d$. Various authors have noted (see \cite{Su08,GuZe01,Be05}) that if $\varrho(\mathsf{A})=1$ then
\[\max_{i_1,\ldots,i_t \in \{0,\ldots,N-1\}} \left\|A_{i_t}\cdots A_{i_1}\right\|\leq Ct^{L_{\mathsf{A}}-1}\]
for some constant $C>0$ depending only on $\mathsf{A}$. (In particular if $L_{\mathsf{A}}=1$ then the corresponding system is marginally stable, an observation which goes back at least to \cite{Ba88}.) This has led several authors to ask in what circumstances a \emph{lower} bound of the same form exists (see for example \cite{ChMaSi12,JuPrBl08}). In this note we will show instead that the sequence \eqref{eq:that-sequence} can oscillate dramatically: we show that there exists a pair of matrices $\mathsf{A}=\{A_0,A_1\}\subset M_3(\R)$ satisfying $\varrho(\mathsf{A})=1$ and $L_{\mathsf{A}}=1$ such that \eqref{eq:that-sequence} grows more slowly than any previously-prescribed unbounded sequence $(a_t)_{t=1}^\infty$ along one sequence of times $t$ and more rapidly than any previously-prescribed sub-linear sequence $(b_t)_{t=1}^\infty=o(t)$ along a complementary sequence of times $t$. We will also show that this behaviour is generic in a certain sense among a particular family of examples.

\section{Previous conjectures and main results}
 In \cite{ChMaSi12}, Y. Chitour, P. Mason and M. Sigalotti noted that the bound $O(t^{L_{\mathsf{A}}-1})$ may in some cases be improved to a smaller polynomial bound $O(t^{\ell_{\mathsf{A}}-1})$ where $\ell_{\mathsf{A}}\leq L_{\mathsf{A}}$ is a positive integer which they called the \emph{resonance degree} of $\mathsf{A}$. In the representation \eqref{eq:bananas} of $\mathsf{A}$ the resonance degree depends only on the diagonal blocks $A_k^{(i,i)}$. Chitour, Mason and Sigalotti conjectured that when these diagonal blocks are fixed and the off-diagonal blocks are chosen freely, for an open and dense set of possible choices of off-diagonal blocks $A^{(i,j)}$ there should exist an initial unit vector $v_0$ and switching law $\sigma$ for which $\|x(t)\| \geq ct^{\ell_{\mathsf{A}}-1}$ for all $t \geq 0$, where $c>0$ does not depend on $t$. Consequently, according to this conjecture the sequence \eqref{eq:that-sequence} should be eventually bounded both above and below by a constant times $t^{\ell_{\mathsf{A}}-1}$ in a similar manner to the situation which results when switching is absent. A similar conjecture was also proposed in the context of continuous-time linear switched systems. In a related earlier work, R.M. Jungers, V.Yu. Protasov and V.D. Blondel similarly asked whether the limit 
\begin{equation}\label{eq:juper}\lim_{t \to \infty} \frac{1}{\log t} \cdot \log \left( \max_{i_1,\ldots,i_t \in \{0,N-1\}} \|A_{i_t}\cdots A_{i_1}\|\right) \end{equation}
always exists when $\mathsf{A}=\{A_0,\ldots,A_{N-1}\}\subset M_d(\R)$ is finite and satisfies $\varrho(\mathsf{A})=1$ (see \cite[Problem 2]{JuPrBl08}).

In the subsequent work \cite{PrJu15} Protasov and Jungers showed that in some cases a discrete-time linear switching system may have resonance degree $2$ yet be marginally stable for all choices of off-diagonal blocks in the representation \eqref{eq:bananas}, demonstrating that the conjecture of Chitour, Mason and Sigalotti is false. 
They also exhibited a pair of matrices $\mathsf{A}=\{A_0,A_1\}\subset M_3(\R)$ such that for a certain constant $C>0$
\[ \max_{i_1,\ldots,i_t \in \{0,1\}} \|A_{i_t}\cdots A_{i_1}\| \leq Ct^{\frac{1}{3}}\]
for all $t\geq 1$, and such that there also exists a trajectory $(x(t))_{t=0}^\infty$ satisfying $\limsup_{t\to \infty} t^{-1/3}\|x(t)\|>0$. 
 (Related examples with an \emph{infinite} set $\mathsf{A}$ had been constructed earlier in \cite{GuZe01}.) These two examples indicated that Chitour, Mason and Sigalotti's bound $O(t^{\ell_{\mathsf{A}}-1})$ can fail to be sharp (moreover, in a manner which is robust with respect to the choice of off-diagonal blocks) but left the existence or nonexistence of the limit \eqref{eq:juper} in general an open question. 
  
In this note we show that the limit \eqref{eq:juper} can entirely fail to exist and that this phenomenon can also occur robustly with respect to changes in the off-diagonal blocks, resolving the question of Jungers, Protasov and Blondel and also providing a new family of counterexamples to the conjecture of Chitour, Mason and Sigalotti in discrete time. Our result is as follows:
\begin{theorem}\label{th:main}
Let $(a_t)_{t=1}^\infty$, $(b_t)_{t=1}^\infty$ be arbitrary sequences of positive real numbers which diverge to infinity and which also satisfy $\lim_{t \to \infty} a_t/t = \lim_{t \to \infty} b_t/t=0$. Let $\lambda \in (-1,1)$ be arbitrary. Then there exists a dense $G_\delta$ set $\Omega \subset [0,2\pi]$ with the following property:  for every $\theta \in \Omega$ there exists a vector space $W\subseteq \R^4$ of dimension at most three such that for every $(u_1,u_2,v_1,v_2) \in\R^4$ the pair of matrices $\mathsf{A}=\{A_0,A_1\}$ defined by
\[A_0:=\begin{pmatrix} \lambda & 0 & u_1 \\
0&-1&u_2\\
0&0&1\end{pmatrix}, \qquad A_1=\begin{pmatrix} \cos \theta &-\sin \theta &v_1 \\ \sin \theta & \cos\theta & v_2 \\ 0&0&1\end{pmatrix} \]
satisfies $\varrho(\mathsf{A})=1$ and satisfies the two properties
\begin{equation}\label{eq:one}\lim_{t \to \infty} \frac{1}{t} \max_{i_1,\ldots,i_t \in \{0,1\}} \left\|A_{i_t}\cdots A_{i_1}\right\|=0,\end{equation}
\begin{equation}\label{eq:two}\liminf_{t \to \infty} \frac{1}{a_t} \max_{i_1,\ldots,i_{t} \in \{0,1\}} \left\|A_{i_{t}}\cdots A_{i_1}\right\|=0,\end{equation}
and such that if $(u_1,u_2,v_1,v_2) \in \R^4 \setminus W$ then additionally
\begin{equation}\label{eq:three}\limsup_{t \to \infty} \frac{1}{b_t} \max_{i_1,\ldots,i_{t} \in \{0,1\}} \left\|A_{i_{t}}\cdots A_{i_1}\right\|=\infty.\end{equation}
Furthermore, the subsequence of times $t$ along which \eqref{eq:two} converges to zero depends only on $\theta$ and not on $(u_1,u_2,v_1,v_2)$.
\end{theorem}
Here we recall that a $G_\delta$ subset of a metric space is by definition any set which is equal to the intersection of a countable family of open sets. A set which admits a dense $G_\delta$ set as a subset is often also called a \emph{residual set} or a \emph{set of second Baire category}. By Baire's theorem the intersection of countably many dense $G_\delta$ subsets of a complete metric space is also a dense $G_\delta$ subset of that space, and this fact is fundamental to the proof of Theorem \ref{th:main}.

  Equations \eqref{eq:two} and \eqref{eq:three} demonstrate that the sequence \eqref{eq:that-sequence} may at once grow arbitrarily slowly along one sequence of times and yet also faster than any prescribed sub-linear sequence along a complementary sequence of times. Thus for generic choices of the parameter $\theta$ there is not only no well-defined rate of polynomial growth for the sequence \eqref{eq:that-sequence}, but this phenomenon moreover holds for an open and dense set of choices of the off-diagonal blocks, contrary to the conjecture of Chitour, Mason and Sigalotti.
  
Choosing $a_t:=1+\log t$ and $b_t:=t/(1+\log t)$ in Theorem \ref{th:main} results in an example where the sequence in \eqref{eq:juper} 
has limit superior equal to $1$ and limit inferior equal to $0$. Since the difference between successive terms in this sequence is $O(1/\log t)$ this implies that rather than having a unique accumulation point as was suggested in \cite{JuPrBl08}, the set of accumulation points of the sequence \eqref{eq:juper} is the entire interval $[0,1]$. The scenario suggested by Jungers, Protasov and Blondel in \cite[Problem 2]{JuPrBl08} therefore fails to hold in what is arguably the strongest possible sense.

We will obtain Theorem \ref{th:main} from the combination of Baire's theorem with the following result:
\begin{theorem}\label{th:easy}
Let $a,b,\theta,\phi \in \R$, $r>0$ and $\lambda \in (-1,1)$. Define two matrices $A_0, A_1 \in M_3(\R)$ by
\[A_0:=\begin{pmatrix}\lambda&0 &a \\0&-1&b\\ 0 &0&1\end{pmatrix},\qquad A_1:=\begin{pmatrix}\cos\theta & -\sin\theta &r \cos\phi \\ \sin\theta &\cos\theta&r\sin\phi\\
0&0 &1\end{pmatrix}\]
and define $\mathsf{A}:=\{A_0,A_1\}$. Then $\varrho(\mathsf{A})=1$ and:
\begin{enumerate}[(i)]
\item\label{it:aly}
The limit
\[\lim_{t \to \infty} \frac{1}{t} \max_{i_1,\ldots,i_t \in \{0,1\}} \left\|A_{i_t}\cdots A_{i_1}\right\| \in [0,\infty)\]
exists.
\item\label{it:hica}
If $\theta=p\pi/q$ with $p$ odd, and if $b\sin\frac{\theta}{2} \neq r\sin\left(\frac{\theta}{2}-\phi\right)$, then $\mathsf{A}$ is marginally unstable and satisfies
\[\lim_{t \to \infty}  \frac{1}{t} \max_{i_1,\ldots,i_t \in \{0,1\}} \left\|A_{i_t}\cdots A_{i_1}\right\|>0.\]
\item\label{it:chy}
If $\theta=p\pi/q$ where $p$ is even, $\gcd(p,q)=1$ and $q>1$, then $\mathsf{A}$ is marginally stable.
\end{enumerate}
\end{theorem}
We remark that the proof of \eqref{it:aly} is in fact very general: it relies only on the fact that the matrices are block upper triangular with exactly two diagonal blocks, and the fact that those diagonal blocks all have norm at most $1$. For reasons of economy we forgo the presentation of a more general statement.

Clauses \eqref{it:hica} and \eqref{it:chy} of Theorem \ref{th:easy} together imply that this family of pairs of matrices exhibits marginal stability on one dense subset of the parameter space, and marginal instability with linear growth of trajectories on a second dense subset of the parameter space. Though earlier results such as \cite{BlTs00} suggest that marginal stability and marginal instability may be challenging to distinguish in practice, we are not aware of any earlier examples to the effect that these two phenomena can both simultaneously be dense in the same region of the parameter space. 

Together with Baire's theorem this dense intermingling of marginal stability and marginal instability with linear growth is essentially the only ingredient necessary to establish the conclusions \eqref{eq:two} and \eqref{eq:three} of Theorem \ref{th:main}. We therefore suggest that the behaviours described in Theorem \ref{th:main} may be quite widespread and should be found in all situations where this intermingling takes place. We in particular anticipate that an analogous phenomenon should exist in continuous-time linear switched systems. These points are discussed further in the conclusions, which are given after the proofs of Theorems \ref{th:easy} and \ref{th:main} below.


\section{Proof of Theorem \ref{th:easy}}
We first prove \eqref{it:aly}. In this part of the proof it will be convenient to write
\[A_0=\begin{pmatrix} B_0& u_0\\ 0&1\end{pmatrix},\qquad A_1=\begin{pmatrix} B_1& u_1\\ 0&1\end{pmatrix} \]
where
\[B_0=\begin{pmatrix} \lambda & 0 \\ 0&-1\end{pmatrix},\qquad u_0=\begin{pmatrix}a\\b\end{pmatrix},\qquad
B_1=\begin{pmatrix} \cos\theta & -\sin\theta\\\sin\theta&\cos\theta\end{pmatrix},\qquad u_1=\begin{pmatrix}r\cos\phi\\r\sin\phi\end{pmatrix}.\]
We note that $\|B_0\|=\|B_1\|=1$ and therefore every possible product $B_{i_t}\cdots B_{i_1}$ has norm at most $1$.  A simple induction shows that for every $t \geq 1$ we have
\[A_{i_t}\cdots A_{i_1} = \begin{pmatrix} B_{i_t}\cdots B_{i_1} & \sum_{k=1}^{t} B_{i_t}\cdots B_{i_{k+1}} u_{i_k}\\ 0&1\end{pmatrix}\]
for all $i_1,\ldots,i_t \in \{0,1\}$. (Here and in the sequel, the empty product $B_{i_t}\cdots B_{i_{t+1}}$ corresponding to the case $k=t$ is understood to be the identity matrix.) For every $t\geq 1$ define
\[\alpha_t:=\max_{i_1,\ldots,i_t \in \{0,1\}} \left\|A_{i_t}\cdots A_{i_1}\right\|,\qquad\beta_t:=\max_{i_1,\ldots,i_t \in\{0,1\}}\left\|\sum_{k=1}^{t} B_{i_t}\cdots B_{i_{k+1}} u_{i_k}\right\|.\]
We will first show that $\lim_{t \to \infty} \frac{1}{t}\beta_t$ exists and is non-negative, and then show that the limit $\lim_{t \to \infty} \frac{1}{t}\alpha_t$ exists and coincides with the former limit.
If $t, \tau \geq 1$ and $i_1,\ldots,i_{t+\tau} \in \{0,1\}$ then
\begin{align*}{\lefteqn{\left\|\sum_{k=1}^{t+\tau} B_{i_{t+\tau}}\cdots B_{i_{k+1}} u_{i_k}\right\| }}&\\
&\leq\left\|\sum_{k=t+1}^{t+\tau} B_{i_{t+\tau}}\cdots B_{i_{k+1}} u_{i_k}\right\|+\left\|\sum_{k=1}^{t} B_{i_{t+\tau}}\cdots B_{i_{k+1}} u_{i_k}\right\|   \\
&=\left\|\sum_{k=1}^{\tau} B_{i_{t+\tau}}\cdots B_{i_{t+k+1}} u_{i_{t+k}}\right\| +\left\|B_{i_{t+\tau}}\cdots B_{i_{t+1}}\cdot \sum_{k=1}^{t} B_{i_{t}}\cdots B_{i_{k+1}} u_{i_k}\right\|   \\
&\leq\left\|\sum_{k=1}^{\tau} B_{i_{t+\tau}}\cdots B_{i_{t+k+1}} u_{i_{t+k}}\right\|+\left\|\sum_{k=1}^{t} B_{i_{t}}\cdots B_{i_{k+1}} u_{i_k}\right\|   \\
&\leq \beta_\tau+\beta_t\end{align*}
where we have used the fact that $\|B_{i_{t+\tau}} \cdots B_{i_{t+1}}\|\leq 1$. Taking the maximum over all possible choices of $i_1,\ldots,i_{t+\tau}$ it follows directly that $\beta_{t+\tau}\leq \beta_\tau+\beta_t$. Since $t,\tau \geq 1$ were arbitrary we have established that $\beta_{t+\tau}\leq \beta_\tau+\beta_t$ for all $t,\tau \geq 1$, and this condition implies the existence of the limit $\lim_{t \to \infty} \frac{1}{t}\beta_t \in [-\infty,\infty)$ by Fekete's subadditivity lemma. The limit is obviously non-negative since the terms $\beta_t/t$ are non-negative. To prove \eqref{it:aly} it is therefore sufficient to show that $\lim_{t \to \infty} \frac{1}{t}|\alpha_t-\beta_t|=0$. If $t\geq 1$ and $i_1,\ldots,i_t \in \{0,1\}$ are arbitrary then we have
\begin{align*} {\lefteqn{\left|\left\|A_{i_t}\cdots A_{i_1}\right\|-\left\|\sum_{k=1}^{t} B_{i_t}\cdots B_{i_{k+1}} u_{i_k}\right\|\right|}}&\\
&=\left|\left\|\begin{pmatrix} B_{i_t}\cdots B_{i_1} & \sum_{k=1}^{t} B_{i_t}\cdots B_{i_{k+1}} u_{i_k} \\ 0&1\end{pmatrix}\right\| - \left\|\begin{pmatrix} 0& \sum_{k=1}^{t} B_{i_t}\cdots B_{i_{k+1}} u_{i_k} \\ 0&0\end{pmatrix}\right\|\right|\\
&\leq \left\|\begin{pmatrix} B_{i_t}\cdots B_{i_1} & \sum_{k=1}^{t} B_{i_t}\cdots B_{i_{k+1}} u_{i_k} \\ 0&1\end{pmatrix} - \begin{pmatrix} 0& \sum_{k=1}^{t} B_{i_t}\cdots B_{i_{k+1}} u_{i_k} \\ 0&0\end{pmatrix} \right\| \\
&=\left\|\begin{pmatrix} B_{i_t}\cdots B_{i_1} & 0 \\ 0&1\end{pmatrix}\right\| \leq 1.\end{align*}
It follows directly that $|\alpha_t-\beta_t|\leq 1$ for all $t \geq 1$ and therefore the desired limit $\lim_{t \to \infty} \frac{1}{t}\alpha_t$ exists and is equal to $\lim_{t \to \infty} \frac{1}{t}\beta_t$. We have proved \eqref{it:aly}.

It follows directly from \eqref{it:aly} that the joint spectral radius of $\mathsf{A}$ cannot be greater than $1$. Since $A_0$ and $A_1$ both clearly have spectral radius $1$, we have
\[\lim_{n \to \infty} \|A_0^t\|^{\frac{1}{t}}=\lim_{t \to \infty} \|A_1^t\|^{\frac{1}{t}}=1\]
by Gelfand's formula. The joint spectral radius of $\mathsf{A}$ thus also cannot be \emph{less} than $1$ and we conclude that $\varrho(\mathsf{A})=1$ as required. We deduce in particular that $\mathsf{A}$ is marginally stable if the products $A_{i_t}\cdots A_{i_1}$ can be bounded independently of $t$ and of the choice of $i_1,\ldots,i_t$, and is marginally unstable otherwise. 

We now prove \eqref{it:hica}. For every $\psi \in \R$ define
\[R_\psi:=\begin{pmatrix} \cos\psi & -\sin \psi \\ \sin\psi&\cos\psi\end{pmatrix},\qquad v_\psi:=\begin{pmatrix}\cos\psi\\\sin\psi\end{pmatrix},\]
so that
\[A_1=\begin{pmatrix} R_\theta & rv_\phi\\ 0&1\end{pmatrix}.\]
A simple induction shows that for every $t \geq 0$
\[A_1^t = \begin{pmatrix} R_{t\theta} & \sum_{k=0}^{t-1}rv_{k\theta+\phi}\\
0&1\end{pmatrix}.\]
Taking $t=q$ yields
\begin{equation}A_1^q = \begin{pmatrix} R_{q\theta} & \sum_{k=0}^{q-1}rv_{k\theta+\phi}\\ 0&1\end{pmatrix}\label{eq:barley}
=\begin{pmatrix} -1&0 &\sum_{k=0}^{q-1}r\cos (k\theta+\phi)\\ 0&-1&\sum_{k=0}^{q-1}r\sin(k\theta+\phi)\\ 0&0&1\end{pmatrix}\end{equation}
since $q\theta=p\pi$ with $p$ odd. Since $p$ is odd and $q$ is an integer, $\theta$ cannot be an integer multiple of $2\pi$. It follows that $e^{i\theta} \neq 1$ and so we may compute
\begin{align*}\sum_{k=0}^{q-1}\cos (k\theta+\phi) & =\Re\left(\sum_{k=0}^{q-1}e^{i(k\theta+\phi)}\right)\\
&=\Re\left(e^{i\phi} \left(\frac{e^{iq\theta}-1}{e^{i\theta}-1} \right)\right)\\
&=\Re\left(e^{i\phi} \left(\frac{e^{ip\pi}-1}{e^{i\theta}-1} \right)\right)\\
&=   \Re\left(-\frac{2e^{i\phi}}{e^{i\theta}-1}\right)\\
&=\Re\left(\left(\frac{2i}{e^{\frac{i\theta}{2}}-e^{\frac{-i\theta}{2}}}\right)\left(\frac{ie^{i\phi}}{e^{\frac{i\theta}{2}}}\right)\right)\\
&=\Re\left(\frac{ie^{i\left(\phi-\frac{\theta}{2}\right)}}{\sin\frac{\theta}{2}}\right)=-\frac{\sin(\phi-\frac{\theta}{2})}{\sin \frac{\theta}{2}}=\frac{\sin(\frac{\theta}{2}-\phi)}{\sin \frac{\theta}{2}}.\end{align*}
Similarly
\[\sum_{k=0}^{q-1} \sin(k\theta+\phi)  =\Im\left(\sum_{k=0}^{q-1}e^{i(k\theta+\phi)}\right)=\Im\left(\frac{ie^{i\left(\phi-\frac{\theta}{2}\right)}}{\sin\frac{\theta}{2}}\right)
=\frac{\cos(\phi-\frac{\theta}{2})}{\sin \frac{\theta}{2}}=\frac{\cos(\frac{\theta}{2}-\phi)}{\sin \frac{\theta}{2}}.\] 
Substituting these values into the formula \eqref{eq:barley} for $A_1^q$ and multiplying on the left by $A_0$ yields
 \[A_0A_1^q = \begin{pmatrix}\lambda&0 &a \\0&-1&b\\ 0 &0&1\end{pmatrix} \begin{pmatrix} -1 &0&\frac{r\cos(\frac{\theta}{2}-\phi)}{\sin \frac{\theta}{2}}\\
0&-1&\frac{r\sin(\frac{\theta}{2}-\phi)}{\sin \frac{\theta}{2}}\\
0&0&1\end{pmatrix}=\begin{pmatrix} -\lambda &0&\frac{\lambda r \cos(\frac{\theta}{2}-\phi)}{\sin \frac{\theta}{2}}+a\\
0&1&-\frac{r\sin(\frac{\theta}{2}-\phi)}{\sin \frac{\theta}{2}}+b\\
0&0&1\end{pmatrix}.\]
Since by hypothesis $b\sin\frac{\theta}{2} \neq r\sin(\frac{\theta}{2}-\phi)$ the lower-right $2\times 2$ block of this matrix is a nontrivial Jordan block. Since additionally $|\lambda|<1$ we therefore obtain
\[\lim_{t \to \infty}\frac{1}{t}\left\|\left(A_0A_1^q\right)^t\right\| = \left|\frac{r\sin(\frac{\theta}{2}-\phi)}{\sin \frac{\theta}{2}}-b\right|>0.\]
It follows using \eqref{it:aly} that
\begin{align*}\lim_{t \to \infty} \frac{1}{t} \max_{i_1,\ldots,i_t \in \{0,1\}} \left\|A_{i_t}\cdots A_{i_1}\right\| &=\lim_{t \to \infty} \frac{1}{t(q+1)} \max_{i_1,\ldots,i_{t(q+1)} \in \{0,1\}} \left\|A_{i_{t(q+1)}}\cdots A_{i_1}\right\|\\
& \geq \lim_{t \to \infty} \frac{1}{t(q+1)} \left\|\left(A_0A_1^q\right)^t\right\|\\
& =\frac{1}{q+1}\left|\frac{r\sin(\frac{\theta}{2}-\phi)}{\sin \frac{\theta}{2}}-b\right| >0\end{align*}
and in particular $\mathsf{A}$ is marginally unstable. We have proved \eqref{it:hica}.

The proof of \eqref{it:chy} is the longest part of the proof and proceeds via a series of short claims. The first such claim is that $A_1^q=I$. Indeed, similarly to \eqref{it:hica} we may calculate
\[A_1^q = \begin{pmatrix} R_{q\theta} & \sum_{k=0}^{q-1}rv_{k\theta+\phi}\\
0&1\end{pmatrix}=\begin{pmatrix} 1&0 &\sum_{k=0}^{q-1}r\cos (k\theta+\phi)\\
0&1&\sum_{k=0}^{q-1}r\sin(k\theta+\phi)\\ 0&0&1\end{pmatrix}\]
since $q\theta$ is now assumed to be an even integer multiple of $\pi$. Since $q>1$ and $\gcd(p,q)=1$ it is impossible for $\theta$ to be an even integer multiple of $\pi$, so $e^{i\theta}\neq 1$. On the other hand $e^{iq\theta}=e^{ip\pi}=1$ since $p$ is even, so we may calculate
\[\sum_{k=0}^{q-1}\cos (k\theta+\phi)  =\Re\left(\sum_{k=0}^{q-1}e^{i(k\theta+\phi)}\right)=\Re\left(e^{i\phi} \left(\frac{e^{iq\theta}-1}{e^{i\theta}-1} \right)\right)=\Re(0)=0\]
and likewise $\sum_{k=0}^{q-1}\sin (k\theta+\phi)=0$. The claimed identity $A_1^q=I$ follows.

Now define
\[P:=\begin{pmatrix}\lambda&0\\0&-1\end{pmatrix},\qquad w:=\begin{pmatrix}a\\b\end{pmatrix}\]
so that
\[A_0=\begin{pmatrix}P&w\\0&1\end{pmatrix}.\]
We observe that $\|P\|=\|R_\theta\|=1$. Our next claim is that there exists $\kappa>0$ such that $\|PR_{\theta}^tP\|\leq 1-\kappa$ for every integer $t \geq 0$ which is not divisible by $q$. Since $R_\theta^q=R_{q\theta}=I$ we may without loss of generality reduce any such integer $t$ modulo $q$ so as to satisfy $1 \leq t<q$ without altering the value of $\|PR_\theta^tP\|$. It thus suffices to show that $\max_{1 \leq t<q} \|PR_{\theta}^tP\|<1$. If $1 \leq t<q$, let $v \in \R^2$ be a unit vector. If $v$ does not lie on the vertical axis then $\|Pv\|<1$ and therefore $\|PR_{\theta}^tPv\|\leq \|Pv\|<1$. On the other hand if $v$ \emph{does} lie on the vertical axis then $R_{\theta}^tPv$ does not, because $tp/q$ is not an integer and $\sin t\theta=\sin (tp\pi/q)$ is therefore not zero. (Here the fact that $\gcd(p,q)=1$ ensures that $q$ does not divide $tp$ and therefore $tp\pi/q$ is not an integer multiple of $\pi$.) It follows that in this case $\|PR_{\theta}^tPv\| <\|R_{\theta }^tPv\|= 1$. We conclude that $\|PR_{\theta}^tP\|<1$ for each $t=1,\ldots,q-1$ and the claim follows.

We now claim that  for all $m \geq 0$,
\[\left\|P^{k_{m}} R_\theta^{\ell_{m}}  \cdots P^{k_2}R_\theta^{\ell_2}P^{k_1}R_\theta^{\ell_1} \right\|\leq (1-\kappa)^{\lfloor \frac{m}{2}\rfloor} \]
for every choice of integers $k_1,\ldots,k_{m} \geq 1$ and $\ell_1,\ldots,\ell_m \geq 1$ such that none of $\ell_1,\ldots,\ell_m$ is divisible by $q$. We will prove the claim using separate inductions on odd and even $m$. The case $m=0$ simply asserts the trivial inequality $\|I\|\leq 1$ and the case $m=1$ asserts that $\|P^{k_1}R_\theta^{\ell_1}\|\leq 1$. These are clearly both true.  It is therefore sufficient to show that if the claim holds for some value of $m$ (which may be either odd or even) then it also holds for $m+2$. Indeed, given $k_1,\ldots,k_{m+2} \geq 1$ and $\ell_1,\ldots,\ell_{m+2} \geq 1$ such that none of $\ell_1,\ldots,\ell_{m+2}$ is divisible by $q$, we have
\begin{align*}{\lefteqn{\left\|P^{k_{m+2}} R_\theta^{\ell_{m+2}}  \cdots P^{k_1}R_\theta^{\ell_1}\right\|}}&\\
&\leq \left\|P^{k_{m+2}} R_\theta^{\ell_{m+2}} P^{k_{m+1}}\right\|\cdot \left\|R_\theta^{\ell_{m+1}}\right\| \cdot \left\|P^{k_{m}} R_\theta^{\ell_{m}}  \cdots P^{k_1}R_\theta^{\ell_1}\right\|\\
&\leq \left\|PR_\theta^{\ell_{m+2}} P\right\| (1-\kappa)^{\lfloor \frac{m}{2}\rfloor}\\
& \leq (1-\kappa)^{1+\lfloor \frac{m}{2}\rfloor}=(1-\kappa)^{\lfloor \frac{m+2}{2}\rfloor}\end{align*}
where we have used the induction hypothesis, the immediately preceding claim, and the elementary inequalities $\|P\|\leq 1$ and $\|R_\theta\|\leq 1$. The result thus follows by induction and the claim is proved.

Since $A_0$ and $A_1$ are both diagonalisable over $\mathbb{C}$ the norms $\|A_0^t\|$ and $\|A_1^t\|$ can be bounded independently of $t$. Choose $C_1>1$ such that $\|A_0^t\|\leq C_1$ and $\|A_1^t\|\leq C_1$ for every $t\geq 0$. For every $\ell,k \geq 0$ define a vector $w_{k,\ell} \in \R^2$ by
\[A_0^kA_1^\ell =\begin{pmatrix} P^kR_\theta^\ell & w_{k,\ell} \\ 0&1\end{pmatrix}.\]
Clearly we have $\|w_{k,\ell}\|\leq \|A_0^kA_1^\ell\|\leq C_1^2$ for every $k,\ell \geq 0$. Our final claim is as follows: there exists a constant $C_2>0$ such that for every $m \geq 0$, for every choice of integers $k_1,\ldots,k_{m} \geq 1$ and $\ell_1,\ldots,\ell_m \geq 1$ such that none of $\ell_1,\ldots,\ell_m$ is divisible by $q$, we have
\begin{equation}\label{eq:ymo}\left\|A_0^{k_{m}}A_1^{\ell_m}\cdots A_0^{k_1} A_1^{\ell_1} \right\| \leq C_2.\end{equation}
By a straightforward induction on $m$ we have
\[A_0^{k_m}A_1^{\ell_m}\cdots A_0^{k_1}A_1^{\ell_1}=\begin{pmatrix}P^{k_m}R_\theta^{\ell_m}\cdots P^{k_1}R_\theta^{\ell_1}& \sum_{j=1}^m P^{k_m}R_\theta^{\ell_m}\cdots P^{k_{j+1}}R_\theta^{\ell_{j+1}} w_{k_j,\ell_j}\\ 0&1\end{pmatrix}\]
for all $k_1,\ldots,k_m \geq 1$ and $\ell_1,\ldots,\ell_m \geq 1$, for all $m \geq 0$. Precisely as in the proof of \eqref{it:aly}, by the reverse triangle inequality the difference 
\[\left|\left\|A_0^{k_m}A_1^{\ell_m}\cdots A_0^{k_1}A_1^{\ell_1} \right\|- \left\|\sum_{j=1}^m P^{k_m}R_\theta^{\ell_m}\cdots P^{k_{j+1}}R_\theta^{\ell_{j+1}} w_{k_j,\ell_j}\right\|\right|\]
is bounded by
\[ \left\|\begin{pmatrix}P^{k_m}R_\theta^{\ell_m}\cdots P^{k_1}R_\theta^{\ell_1}&0\\ 0&1\end{pmatrix}\right\|\leq1,\]
so in particular if none of the integers $\ell_i$ is divisible by $q$,
\begin{align*}\left\|A_0^{k_m}A_1^{\ell_m}\cdots A_0^{k_1}A_1^{\ell_1}\right\|
 &\leq 1+\sum_{j=1}^m \left\|P^{k_m}R_\theta^{\ell_m}\cdots P^{k_{j+1}}R_\theta^{\ell_{j+1}} w_{k_j,\ell_j}\right\|\\
& \leq 1+\sum_{j=1}^m C_1^2\left\|P^{k_m}R_\theta^{\ell_m}\cdots P^{k_{j+1}}R_\theta^{\ell_{j+1}}\right\|\\
& \leq 1+\sum_{j=1}^m C_1^2\left(1-\kappa\right)^{\left\lfloor \frac{m-j}{2}\right\rfloor}\\
& = 1+\sum_{j=0}^{m-1} C_1^2\left(1-\kappa\right)^{\left\lfloor \frac{j}{2}\right\rfloor}\\
&<1+\sum_{j=0}^\infty C_1^2\left(1-\kappa\right)^{\left\lfloor \frac{j}{2}\right\rfloor}\\
& = 1+2\sum_{j=0}^\infty C_1^2\left(1-\kappa\right)^j =1+\frac{2C_1^2}{\kappa}\end{align*}
and the claim is proved with $C_2:=1+2C_1^2\kappa^{-1}$.

We may now prove that $\mathsf{A}$ is marginally stable. Consider a product $M:=A_{i_t}\cdots A_{i_1}$ of the matrices $A_0$ and $A_1$.
Write $M$ in the form
\[M=A_1^{\ell_{m+1}}A_0^{k_{m}}A_1^{\ell_m} \cdots A_0^{k_1}A_1^{\ell_1}A_0^{k_0} \]
where $k_1,\ldots,k_{m}, \ell_1,\ldots,\ell_m \geq1$, $k_0, \ell_{m+1}\geq 0$ and where $m \geq 0$ is chosen as small as possible. If $m\geq 1$ then we observe that no integer $\ell_j$ with $j\in \{1,\ldots,m\}$ can be divisible by $q$, since if this were the case then we would have $A_0^{k_{i}}A^{\ell_j}_1A_0^{k_{j-1}}=A_0^{k_{j}+k_{j-1}}$ by virtue of the identity $A_1^q=I$ and it follows that there would exist a representation of $M$ in the above form with smaller $m$, contradicting the minimality of $m$. The previous claim is therefore applicable to the product $A_0^{k_m}A_1^{\ell_m}\cdots A_0^{k_1}A_1^{\ell_1}$ and we may apply \eqref{eq:ymo} to deduce that
\[\|M\| \leq C_1^2 \left\|A_0^{k_m} A_1^{\ell_m}\cdots A_0^{k_1}A_1^{\ell_1}\right\| \leq C_1^2C_2.\]
Since $t \geq 1$ and $i_1,\ldots,i_t$ were arbitrary this shows that $\mathsf{A}$ is marginally stable. The proof of the Theorem is complete.

\section{Proof of Theorem \ref{th:main}}
In this section it will be convenient for our notation to identify $\R^4$ with $\R^2\times\R^2$ in the obvious manner. For every $\theta \in \R$ and $u,v \in \R^2$ define
\[A_0^{\theta,u,v}:=\begin{pmatrix} P&u\\0&1\end{pmatrix},\qquad A_1^{\theta,u,v}:=\begin{pmatrix} R_\theta & v\\0&1\end{pmatrix}\]
where the notation $P$ and $R_\theta$ is as introduced in the previous section. Let $e_1,e_2$ denote the standard basis vectors in $\R^2$. Let $w, w' \in \R^2$ be arbitrary vectors with the second co-ordinate of $w'$ being nonzero. Define $\Omega_1\subseteq [0,2\pi]$ to be the set of all $\theta \in [0,2\pi]$ such that
\[ \liminf_{t \to \infty} \frac{1}{a_t}  \max_{i_1,\ldots,i_t \in \{0,1\}} \max_{j,k\in \{1,2\}}\left\|A_{i_t}^{\theta,e_j,e_k}\cdots A^{\theta,e_j,e_k}_{i_1}\right\|=0\]
and let  $\Omega_2\subseteq [0,2\pi]$ denote the set of all $\theta \in [0,2\pi]$ such that
\[ \limsup_{t \to \infty} \frac{1}{b_t} \max_{i_1,\ldots,i_t\in\{0,1\}} \left\|A_{i_t}^{\theta,w,w'}\cdots A^{\theta,w,w'}_{i_1}\right\|=\infty.\]
We will show that $\Omega:=\Omega_1 \cap \Omega_2$ has the properties required in the Theorem. 

We first show that each of $\Omega_1$ and $\Omega_2$ is a $G_\delta$ set.  For each $t,r \in \N$ let $U_{t,r}\subset [0,2\pi]$ denote the set
\[\left\{\theta \in [0,2\pi] \colon \frac{1}{a_t} \max_{i_1,\ldots,i_t}  \max_{j,k\in\{1,2\}} \left\|A_{i_t}^{\theta,e_j,e_k}\cdots A^{\theta,e_j,e_k}_{i_1}\right\|<\frac{1}{r}\right\}\]
and let $V_{n,r}\subset [0,2\pi]$ denote the set
\[ \left\{\theta \in [0,2\pi]\colon \frac{1}{b_t} \max_{i_1,\ldots,i_t \in \{0,1\}} \left\|A_{i_t}^{\theta,w,w'}\cdots A^{\theta,w,w'}_{i_1}\right\|>r\right\}.\]
For each fixed $t \geq 1$ the functions
\[\theta \mapsto \max_{i_1,\ldots,i_t \in \{0,1\}} \max_{j,k\in\{1,2\}} \left\|A_{i_t}^{\theta,e_j,e_k}\cdots A^{\theta,e_j,e_k}_{i_1}\right\|\]
and
\[\theta \mapsto \max_{i_1,\ldots,i_t \in \{0,1\}}  \left\|A_{i_t}^{\theta,w,w'}\cdots A^{\theta,w,w'}_{i_1}\right\|\]
are clearly continuous, so each $U_{t,r}$ is an open subset of $[0,2\pi]$ and so too is every $V_{t,r}$. Now we have $\theta \in \Omega_1$ if and only if for every $r \geq 1$ we have $\theta \in U_{t,r}$ for infinitely many $t$, which is to say
\[\Omega_1=\bigcap_{r=1}^\infty \bigcap_{m=1}^\infty \bigcup_{t=m}^\infty U_{t,r}.\]
Likewise
\[\Omega_2=\bigcap_{r=1}^\infty \bigcap_{m=1}^\infty \bigcup_{t=m}^\infty V_{t,r}.\]
It follows that each of $\Omega_1$ and $\Omega_2$ is a $G_\delta$ set as required.

We next show that each of $\Omega_1$ and $\Omega_2$ is dense. Since the sequence $(a_t)$ is unbounded, the set of all $\theta \in [0,2\pi]$ such that $\{A_0^{\theta, e_j, e_k}, A_1^{\theta, e_j,e_k}\}$ is marginally stable for all choices of $j,k\in\{1, 2\}$ is a subset of $\Omega_1$. By Theorem \ref{th:easy}\eqref{it:chy} that set includes the set
\[\left\{\frac{p\pi}{q} \in (0,2\pi) \colon p\text{ is even, }q>1\text{ is odd and }\gcd(p,q)=1\right\}\]
which is dense in $[0,2\pi]$, and it follows that $\Omega_1$ is a dense subset of $[0,2\pi]$. Let us write
\[w=\begin{pmatrix}\alpha\\\beta\end{pmatrix},\qquad w'=\begin{pmatrix} r\cos \phi\\ r\sin\phi\end{pmatrix}\]
where $r>0$. Since $b_t=o(t)$ the set of all $\theta \in [0,2\pi]$ such that
\[\lim_{t \to \infty} \frac{1}{t}\max_{i_1,\ldots,i_t\in \{0,1\}} \left\|A_{i_t}^{\theta,w,w'} \cdots A_{i_1}^{\theta,w,w'}\right\|>0 \]
is a subset of $\Omega_2$, and  by Theorem \ref{th:easy}\eqref{it:hica} this set in turn contains every $\theta=p\pi/q \in (0,2\pi)$ such that $p$ is odd and such that $\beta\sin \left(\frac{\theta}{2}\right) \neq r\sin\left(\frac{\theta}{2}-\phi\right)$. Since by hypothesis the second co-ordinate of $w'$ is nonzero, we have $\sin \phi \neq 0$. By an easy calculus exercise it follows that the function $\theta \mapsto \sin(\frac{\theta}{2}-\phi) /\sin(\frac{\theta}{2})$ is strictly monotone on $(0,2\pi)$ and therefore the equation $\beta\sin \left(\frac{\theta}{2}\right) = r\sin\left(\frac{\theta}{2}-\phi\right)$ can hold for at most one $\theta \in (0,2\pi)$. In particular the set of all  $\theta=p\pi/q \in (0,2\pi)$ such that $p$ is odd and such that $\beta\sin \left(\frac{\theta}{2}\right) \neq r\sin\left(\frac{\theta}{2}-\phi\right)$ is a dense subset of $[0,2\pi]$. We conclude that $\Omega_2$ is also a dense subset of $[0,2\pi]$. 

We have now shown that $\Omega_1$ and $\Omega_2$ are dense $G_\delta$ subsets of $[0,2\pi]$, and therefore by Baire's theorem the set $\Omega:=\Omega_1 \cap\Omega_2$ is also a dense $G_\delta$ subset of $[0,2\pi]$. It remains to show that every $\theta \in \Omega$ has the properties required by the Theorem. Fix $\theta \in \Omega$ for the remainder of the proof. Since by definition we have
\[ \liminf_{t \to \infty} \frac{1}{a_t}  \max_{i_1,\ldots,i_t \in \{0,1\}} \max_{j,k=1,2}\left\|A_{i_t}^{\theta,e_j,e_k}\cdots A^{\theta,e_j,e_k}_{i_1}\right\|=0\]
it follows that there exists a sequence of natural numbers $(t_n)_{n=1}^\infty$ diverging to infinity such that 
\[ \lim_{n \to \infty} \frac{1}{a_{t_n}}  \max_{i_1,\ldots,i_{t_n} \in \{0,1\}} \max_{j,k=1,2}\left\|A_{i_{t_n}}^{\theta,e_j,e_k}\cdots A^{\theta,e_j,e_k}_{i_1}\right\|=0.\]
Let $W_1$ be the set of all $(u,v) \in \R^2 \times \R^2 \simeq \R^4$ such that 
\[ \limsup_{n \to \infty} \frac{1}{a_{t_n}}  \max_{i_1,\ldots,i_{t_n} \in \{0,1\}}\left\|A_{i_{t_n}}^{\theta,u,v}\cdots A^{\theta,u,v}_{i_1}\right\|=0\]
and $W_2$ the set of all $(u,v) \in \R^2 \times \R^2 \simeq \R^4$ such that 
\[ \limsup_{t \to \infty} \frac{1}{b_t} \max_{i_1,\ldots,i_t\in\{0,1\}} \left\|A_{i_t}^{\theta,u,v}\cdots A^{\theta,u,v}_{i_1}\right\|<\infty.\]
We claim that $W_1$ and $W_2$ are both vector spaces. Given two vectors $w_0, w_1 \in \R^2$ write
\[A_0^{\theta, w_0, w_1} = \begin{pmatrix} B_0 & w_0 \\0 &1\end{pmatrix}, \qquad A_1^{\theta, w_0, w_1} = \begin{pmatrix} B_1 & w_1 \\0 &1\end{pmatrix}\]
where $B_0:=P$, $B_1:=R_\theta$. Similarly to the proof of Theorem \ref{th:easy}, for every fixed $t \geq 1$ and $i_1,\ldots,i_t \in \{0,1\}$ we have
\[A_{i_t}^{\theta,w_0,w_1}\cdots A_{i_1}^{\theta,w_0,w_1} = \begin{pmatrix} B_{i_t}\cdots B_{i_1} & \sum_{k=1}^{t} B_{i_t}\cdots B_{i_{k+1}} w_{i_k}\\ 0&1\end{pmatrix}.\]
Since the diagonal terms are bounded independently of $t$ and $i_1,\ldots,i_t$, and since $(a_t)$ and $(b_t)$ diverge to infinity, it follows directly that the values
\[\limsup_{n \to \infty} \frac{1}{a_{t_n}}  \max_{i_1,\ldots,i_{t_n} \in \{0,1\}}\left\|A_{i_{t_n}}^{\theta,w_0,w_1}\cdots A^{\theta,w_0,w_1}_{i_1}\right\| \in [0,\infty]\]
and
\[\limsup_{t \to \infty} \frac{1}{b_t} \max_{i_1,\ldots,i_t\in\{0,1\}} \left\|A_{i_t}^{\theta,w_0,w_1}\cdots A^{\theta,w_0,w_1}_{i_1}\right\|\in [0,\infty]\]
coincide respectively with
\[\limsup_{n \to \infty} \frac{1}{a_{t_n}}  \max_{i_1,\ldots,i_{t_n} \in \{0,1\}}\left\| \sum_{k=1}^{t_n} B_{i_{t_n}}\cdots B_{i_{k+1}} w_{i_k}\right\|\]
and
\[\limsup_{t \to \infty} \frac{1}{b_t} \max_{i_1,\ldots,i_t\in\{0,1\}} \left\| \sum_{k=1}^{t} B_{i_t}\cdots B_{i_{k+1}} w_{i_k}\right\|.\]
These latter are both easily seen to be convex and homogenous functions of the variable $(w_0,w_1) \in \R^2 \times \R^2 \simeq \R^4$. The set on which a convex and homogenous function $\R^4 \to [0,\infty]$ is finite, and the set on which such a function is zero, are obviously vector spaces. Therefore $W_1$ and $W_2$ are vector subspaces of $\R^4$ as required. By construction $W_1$ contains each of the four pairs $(e_j,e_k)$ for $j,k=1,2$ and $W_2$ does \emph{not} contain the pair $(w,w')$. Consequently $W_1$ is four-dimensional and hence is equal to $\R^4$, while $W_2$ is at most three-dimensional. 
 
 We may now complete the proof. If $(u_1,u_2,v_1,v_2) \in \R^4$ is arbitrary then it belongs to $W_1=\R^4$ and therefore \eqref{eq:two} is satisfied by the definition of $W_1$. By Theorem \ref{th:easy}\eqref{it:aly} the limit in \eqref{eq:one} exists, and since $a_t=o(t)$ it follows from \eqref{eq:two} that this limit must be zero.  If $(u_1,u_2,v_1,v_2) \notin W_2=:W$ then \eqref{eq:three} holds by the definition of $W_2$. The proof of the Theorem is complete.

 \section{Conclusions}
 
We have constructed examples of three-dimensional discrete-time linear switched systems with two switching states for which the maximal uniform rate of growth of unstable trajectories is highly irregular, alternating between arbitrarily slow growth at some times and arbitrarily fast sub-linear growth at other times. This suggests to us that there can be no hope of determining the uniform growth rate of a marginally unstable discrete-time linear switched system simply by observing the variety of its behaviour up to some finite time threshold.
 We have also shown that such examples are robust with respect to the choice of their off-diagonal matrix blocks, contrary to a conjecture of Chitour, Mason and Sigalotti \cite[Conjecture 1]{ChMaSi12}. This demonstrates that the sequence \eqref{eq:that-sequence} does not either typically or robustly grow at a precise polynomial rate in the manner hoped for in \cite{ChMaSi12,JuPrBl08}. It remains unclear whether any sharper upper bound for \eqref{eq:that-sequence} than that given in \cite{ChMaSi12} might hold in any significant generality.
 
 Significantly, the existence of discrete linear switched systems with the properties described in Theorem \ref{th:main} arose as a  direct mathematical consequence of the fact that marginally stable systems and marginally unstable systems with linear growth of trajectories were found to be both simultaneously dense in the same region of parameter space. We therefore anticipate that marginally unstable systems with this kind of irregular trajectory growth should occur densely in \emph{any} situation in which both marginally stable systems and linearly marginally unstable systems coexist densely in the same parameter region.
We in particular anticipate that the same phenomenon should occur for continuous-time linear switched systems:
\begin{conjecture}\label{co:xs}
There exist integers $k,d \geq 1$ such that given any two increasing  functions $\alpha,\beta \colon [0,\infty) \to (0,\infty)$ which diverge to infinity and satisfy $\lim_{t\to\infty} \alpha(t)/t=\lim_{t \to \infty} \beta(t)/t=0$ we may choose matrices $A_1,\ldots,A_k \in M_d(\R)$ with the following property. For every unit vector $u \in \R^d$ and every measurable function $L$ from $[0,\infty)$ to the convex hull of $\{A_1,\ldots,A_k\}$, let $v_{L,u} \colon [0,\infty) \to \R^d$ be the unique Lipschitz continuous function satisfying $v_{L,u}'(t)=L(t)v(t)$ a.e. and $v_{L,u}(0)=u$. Then
\[\lim_{t \to \infty} \frac{1}{t} \sup_{L}\sup_{u} \|v_{L,u}(t)\|=0,\]
\[\liminf_{t \to \infty} \frac{1}{\alpha(t)} \sup_{L}\sup_{u} \|v_{L,u}(t)\|=0\]
and
\[\limsup_{t \to \infty} \frac{1}{\beta(t)} \sup_{L}\sup_{u} \|v_{L,u}(t)\|=\infty.\]
Moreover these properties persist for an open and dense set of choices of off-diagonal blocks in a maximal block upper triangular representation of $A_1,\ldots,A_k$.
\end{conjecture}
If a family of continuous-time linear switched systems could be found which depended continuously on a parameter in such a way that both marginal stability and marginal instability with linear rate of trajectory growth were both dense in the parameter space, Conjecture \ref{co:xs} would then follow by arguments broadly analogous to those in this article. In particular, standard results on the weak-* continuity of solutions as a function of the switching law (such as are found in \cite{LiSu99,CoKl00}) should be sufficient to prove that the desired set of parameters is a $G_\delta$ set. The difficulty in proving Conjecture \ref{co:xs} would therefore seem to lie mainly in constructing a family of examples for which marginal stability and instability are densely intermingled in the parameter space. In the discrete-time case this was achieved via resonances between the rotation parameter and the lengths of the intervals between the discrete values of the time parameter $t=1,2,3,\ldots$ itself. In the continuous-time case this part of the construction seems more problematic since the latter intervals no longer exist, but a resonance serving the same function might perhaps be achieved by pairing together two systems with the property that the most-unstable switching law is constrained to be a bang-bang periodic control (such as in Case CC.2.2 of \cite[Theorem 2.3]{Bo02}) in such a way that the periods of the two most-unstable laws resonate with one another.  

The results in this article do not exhaustively describe the possible behaviours of the system considered in Theorem \ref{th:easy}. It is possible to show that 
\[\lim_{t \to \infty} \frac{1}{t}\max_{i_1,\ldots,i_t\in \{0,1\}} \|A_{i_t}\cdots A_{i_1}\|=0\]
whenever the parameter $\theta$ is an irrational multiple of $\pi$, but we are only able to prove this using technical tools from ergodic theory (specifically, those in the appendix to \cite{Mo13}) which the present article otherwise does not require, and since this result has no relevance in establishing Theorem \ref{th:main} it is omitted. 

We also have not investigated the question of whether $\mathsf{A}$ is marginally stable when $\theta=p\pi/q$ with $p$ odd but when the condition $r\sin\frac{\theta}{2} \neq b\sin\left(\phi-\frac{\theta}{2}\right)$ is not satisfied. We believe that marginal stability of $\mathsf{A}$ in this latter case is likely. More speculatively we ask the following two questions:
\begin{question}\label{qu:ip}
For the system $\mathsf{A}$ defined in Theorem \ref{th:easy}, what can be said about the behaviour of the sequence
\begin{equation}\label{eq:ym}\max_{i_1,\ldots,i_t \in \{0,1\}} \|A_{i_t}\cdots A_{i_1}\|\end{equation}
for Lebesgue almost every value of $\theta \in (0,2\pi)$?
\end{question}
\begin{question}\label{qu:ackslikeaduck}
Let $(a_t)_{t=1}^\infty$ be an arbitrary increasing and unbounded sequence of positive real numbers. Do there exist values of the parameters $\theta$, $\lambda$, $a$, $b$, $r$, $\phi$ in Theorem \ref{th:easy} such that $\mathsf{A}$ is marginally unstable, but also satisfies
\[\limsup_{t\to \infty} \frac{1}{a_t}\max_{i_1,\ldots,i_t\in\{0,1\}}\|A_{i_t}\cdots A_{i_1}\|=0?\]
\end{question}
In general, behaviours which occur for generic parameters in the sense of Baire category can contrast strongly with the behaviours which occur for generic parameters in the sense of measure theory; examples of this contrast drawn from a range of topics across mathematical analysis may be found in the monograph \cite{Ox71}. Closer to the topic of this article, continuity of the lower spectral radius occurs on a set of second Baire category but is conjectured to fail on a set of positive Lebesgue measure: see \cite[\S7.3]{BoMo15}. Regarding Question \ref{qu:ip} we therefore suggest that the behaviour of the sequence \eqref{eq:ym} for values of $\theta$ which are typical in the sense of Lebesgue measure might contrast strongly with the situation of Theorem \ref{th:main}, and for Lebesgue-typical values of $\theta$ the example in Theorem \ref{th:main} might obey a simple polynomial growth pattern such as $\sqrt{t}$ or $t^{1/3}$. Question \ref{qu:ackslikeaduck} on the other hand is a specialisation of an earlier question of Protasov and Jungers \cite[\S5]{PrJu15}.

We lastly note that in this article we have addressed only the growth of $\|A_{i_t}\cdots A_{i_1}\|$ uniformly with respect to the choice of switching law and initial vector $v_0$. It is natural to ask whether Theorem \ref{th:main} can be modified so as to refer to pointwise growth with respect to the choice of switching law and initial vector, in those parts where this produces a stronger statement:
\begin{question}
Let $(a_t)_{t=1}^\infty$, $(b_t)_{t=1}^\infty$ be arbitrary unbounded sequences of positive real numbers such that $\lim_{t \to \infty} a_t/t = \lim_{t \to \infty} b_t/t=0$. Do there exist $d \geq 1$, $N \geq 2$, a finite marginally unstable set of matrices $\mathsf{A}=\{A_0,\ldots, A_{N-1}\}\subset M_d(\R)$, a unit vector $v_0 \in \R^d$ and a switching sequence $(z_t)_{t=1}^\infty \in \{0,1,\ldots,N-1\}^\N$ such that
\[\limsup_{t \to \infty}\frac{1}{b_t}\|A_{z_t}\cdots A_{z_1}v_0\|=\infty\]
and also
\[\liminf_{t \to \infty}\frac{1}{a_t} \max_{i_1,\ldots,i_t \in \{0,1,\ldots,N-1\}} \|A_{i_t}\cdots A_{i_1}\|=0,\]
\[\lim_{t \to \infty}\frac{1}{t}\max_{i_1,\ldots,i_t \in \{0,1,\ldots,N-1\}} \|A_{i_t}\cdots A_{i_1}\|=0?\]
\end{question}

 \section{Acknowledgements}
 
This research was partially supported by the Leverhulme Trust (Research Project Grant RPG-2016-194).
 
\bibliographystyle{acm}
\bibliography{swatch}

\end{document}